\documentclass[12pt]{article}
\usepackage{amsmath, amsfonts, amssymb}

\newtheorem{theorem}{Theorem}

\newtheorem{lemma}[theorem]{Lemma}
\def\Proof{\noindent{\it Proof.}\ }

\def\N{{\mathbb N}}

\def\C{{\mathbb C}}

\def\l{\lambda}

\def\Sym{\mathfrak{Sym}}
\def\mfS{{\mathfrak S}}
\def\x{{\mathbf x}}
\def\y{{\mathbf y}}
\def\z{{\mathbf z}}

\def\tF{\widetilde{F}}

\def\carre{\Cup}
\def\d{\partial}
\def\moins{\raise 1pt\hbox{{$\scriptstyle -$}}}
\def\plus{\raise 1pt\hbox{{$\scriptstyle +$}} }

\def\bmoins{\fbox{\raise 1pt\hbox{{$\scriptstyle -$}}1}}

\def\QED{\hfill Q.E.D.}

\newdimen\unit
\def\o{$\scriptscriptstyle{{\rm o}}$}
\def\grape(#1,#2)#3{\raise#2\unit\rlap{\kern#1\unit #3}\ignorespaces}

\def\gg{{\unit=1mm
\hbox {\grape(2,1.2){'}
       \grape(1,2)\o
       \grape(2,2)\o
       \grape(3,2)\o
       \grape(1.5,1)\o
       \grape(2.5,1)\o
       \grape(2,0)\o\
}\kern 3.5 \unit}}

\def\gfill{\leaders\hbox to 1.2em{\hss\gg\hss}\hfill}
\def\frise{\medskip \centerline{\hbox to 8cm{\gfill} }\bigskip}

\begin{document}

\centerline{\Large\bf Gaudin functions, and Euler-Poincar\'e characteristics}

\medskip
\centerline{\large\emph{ Alain Lascoux}}

\bigskip

\frise
\begin{abstract}
Given two positive integers $n,r$, we define the Gaudin function
of level $r$ to be quotient of the numerator of
$$\det\Bigl(\,  ((x_i-y_j)(x_i-ty_j) \cdots (x_i-t^r y_j))^{-1} 
                                                \Bigr)_{i,j=1\ldots n}$$
by the two Vandermonde in  $x$ and $y$.
We show that it can be characterized by  specializing the $x_i$
into the $y_j$ variables, multiplied by powers of $t$.
This allows us to obtain the Gaudin function of level $1$ (due to Korepin and
Izergin) as the image of a resultant under the 
the Euler-Poincar\'e characteristics of the flag manifold. 
As a corollary, we recover a result of Warnaar about the generating function
of Macdonald polynomials.
\end{abstract}

\frise

\section{Gaudin functions of arbitrary level}

Let $\x=\{ x_1,\ldots, x_n\}$,  $\y=\{ y_1,\ldots, y_n\}$
be two sets of indeterminates of the same cardinality $n$.
The Cauchy determinant $\det\bigl( (x-y)^{-1}  \bigr)_{x\in \x, y\in\y}$
plays a central r\^ole in the theory of symmetric functions \cite{Macdonald}.

Generalizations of this determinant appear in the calculation of
correlation functions of different physical models.
Gaudin \cite[Ch.IV]{Gaudin} obtained the determinant
$$ \det\bigl( (x-y)^{-1} (x-y+\gamma)^{-1}  \bigr)_{x\in \x, y\in\y} \, ,$$
$\gamma$ a parameter, for a bose gas in one dimension. 
Izergin and Korepin \cite{Izergin} solved the 
Heisenberg XXZ-antiferromagnetic model with the help of
$$  \det\bigl( (x-y)^{-1} (x- t y)^{-1}  \bigr)_{x\in \x, y\in\y} \, ,$$
and Kirillov and Smirnov \cite[Th.1]{Kirillov} wrote more general
determinants. These different determinants have led to an abundant literature,
in connection with different statistical models and combinatorial
enumerations, for example the enumeration 
of alternating sign matrices displaying some symmetries\cite{Okada}.

In this text, we shall be interested in a purely algebraic generalization,
having no physical interpretation to offer. 
With $\x$ and $\y$ as above, let $t$ be an extra indeterminate
and $r$ be a positive integer.

We propose to study the determinant 
$$ \det \left(\,\frac{1}{ (x-y)(x-ty)\cdots (x-t^r y)} \, 
                                    \right)_{x\in \x, y\in\y} \ . $$
This is a rational function, which can be written
\begin{equation}
   \frac{\Delta(\x) \Delta(\y)}{R(\x, \y(1+\cdots+t^r))}\,
 F^r_n(\x,\y)  \, ,
\end{equation}
where  $F^r_n(\x,\y)$ is a polynomial symmetrical in both
$x_1,\ldots, x_n$ and $y_1,\ldots, y_n$, 
$\Delta(\x)$ being the Vandermonde $\prod_{i<j} (x_j-x_i)$, the resultant 
$R(\x,\y)$, for two finite sets of indeterminates $\x,\y$, being
the product of differences \\ $ \prod_{x\in\x}\prod_{y\in\y} (x-y)$,
and $\y(1+\cdots+t^r)$ being 
$\{ t^i y_j :\, 1\leq j\leq n,\, 0\leq i\leq r\}$.

\medskip 
We shall call the function $F_n^r(\x,\y)$ the \emph{Gaudin function of 
level $r$}, in recognition of the pioneering work of Gaudin so
well illustrated in his thesis \cite {GaudinThese} and his book \cite{Gaudin}.

Let us introduce $n$ sets of indeterminates
$$  \y^j :=\{ y_j^0=y_j, y_j^1,\ldots, y_j^r \}\ ,\ j=1,\ldots, n \, .$$

Recall that, for any pair $(a,b)$ of indeterminates,
the divided difference relative to $a,b$ is the operator (denoted on the right)
$$   f(a,b) \to  f(a,b)\d_{a,b} = \bigl(f(a,b)-f(b,a) \bigr)\, (a-b)^{-1} \, .$$
 
The determinant $M_n= \bigl| R(x_i,\y^j)^{-1} \bigr|$
is the image of the Cauchy determinant 
under the product of divided differences
$$  \prod_{j=1}^n \d_0^j \cdots \d_{r-1}^j  \, ,$$
where $\d_i^j$ is  relative to the pair
 $ y_j^i, y_j^{i+1}$.  

Indeed, for any $j$, any $x$, one has 
$$ (x-y^0_j)^{-1} \stackrel{\d_0^j}{\longrightarrow} 
(x-y^0_j)^{-1}(x-y^1_j)^{-1} \stackrel{\d_1^j}{\longrightarrow}
(x-y^0_j)^{-1}(x-y^1_j)^{-1}(x-y^2_j)^{-1}  \cdots  \, ,
$$
so that  the divided differences relative to the alphabet
$\y^j$ transform all the terms $(x_i-y^0_j)^{-1}$, inside the Cauchy determinant,
 into $(x_i-y^0_j)^{-1}\cdots (x_i-y^r_j)^{-1}$.

On the other hand, it is classical, and immediate,
 that the Cauchy determinant be equal 
to $\Delta(\x) \Delta(\y)/R(\x,\y)$. It can also be expressed as
$$ \Delta(\x)\, \left|\frac{1}{R(\x,y)},\, \frac{y}{R(\x,y)},  
 \ldots,\,   \frac{y^{n-1}}{R(\x,y)},  \right|_{y=y_1,\ldots,y_n}  \, .$$

Its product by $\prod_i R(\x, \y^i)$ can therefore be written as
\begin{equation}  \label{Cauchy2}
 \Delta(\x)\, \Bigl| (y_i^0)^j\, R(\x, \y^i-y_i^0) \Bigr|_{i=1,\ldots,n,\, 
  j=0,\ldots,n-1}
\, . 
\end{equation}

Since each factor $R(\x,\y^i)$ commutes with any divided difference
$\d_k^j$, $k=0,\ldots, r\moins 1$, $j=1,\ldots,n$, the image of 
(\ref{Cauchy2}) under $\prod_i \d_0^i \cdots \d_{r-1}^i $ 
is equal to $M_n\, \prod_i R(\x, \y^i)$:
\begin{equation}  \label{Cauchy3}
  \frac{M_n}{\Delta(\x) } \, \prod_i R(\x, \y^i) = 
        \Bigl| (y_i^0)^j\, R(\x, \y^i-y_i) \d_0^i \cdots \d_{r-1}^i
  \Bigr|  \, .
\end{equation}

We shall now have recourse to symmetric functions, referring to the last
section for more details, as well as \cite{Cbms,Macdonald}.
For any non negative integer $j$, one has 
$$  (y_i^0)^j\, R(\x, \y^i-y_i^0) \d_0^i \cdots \d_{r-1}^i
  = S_{j;\square}( \y^i ; \y^i -\x) \, , $$
where $\square = [\underbrace{n\moins 1,\ldots,n\moins 1}_r ] = (n\moins 1)^r$,
and $S_{j;\square}$ is a multi-Schur function.

Writing $\y^i = (\y^i -\x) +\x$, and expanding by linearity, then 
$$  S_{j;\square}( \y^i ; \y^i -\x) =
 \sum_{k=0}^{j} S_k(\x) S_{j-k,\, \square}(\y^i -\x) \, .$$
This identity allows us to transform  (\ref{Cauchy3})  into
\begin{equation}  \label{Cauchy4}
 \frac{M_n}{\Delta(\x) } \, \prod_i R(\x, \y^i) =
        \Bigl|  S_{j,\, \square}( \y^i -\x) \Bigr|  \, ,
\end{equation}
where now the entries of row $i$ are Schur functions in the difference 
$\y^i -\x$.

Going back to the original variables, that is, 
specializing each $\y^i$ into $y_i+ty_i+\cdots +t^r y_i$, we obtain
the following expression of the Gaudin function of level $r$.

\begin{theorem}   \label{GaudinMatLevel}
The function 
$$ F^r_n(\x, \y) :=
  \frac{ R(\x, \y(1\plus \cdots \plus t^r))}{\Delta(\x)\Delta(\y)} \,
 \det \left(\,\frac{1}{ (x_i-y_j)(x_i-ty_j)\cdots (x_i-t^r y_j)} \, \right)
$$
is equal to 
$$\frac{1}{\Delta(\y)}\, 
   \Bigl|  S_{j,\, \square}(y_i(1\plus \cdots \plus t^r) -\x) 
      \Bigr|_{i=1,\ldots,n,\, j=0,\ldots,n-1} \, .$$
\end{theorem}

The original case of Izergin, Korepin is for level $1$, and reads
\begin{multline}
 F^1_n(\x, \y)  =  \frac{1}{\Delta(\y)}\, 
   \Bigl|  S_{n-1}(y_i(1\plus t) -\x),\, \\
          S_{1,n-1}(y_i(1\plus t) -\x),\, \ldots,
          S_{n-1,n-1}(y_i(1\plus t) -\x) \,
      \Bigr|_{i=1,\ldots,n} \, .
\end{multline}
For example, for $r=1$, $n=3$, the function is 
\begin{equation*}
 F^1(\x,\y) \Delta(\y)  = \begin{vmatrix}
S_{022}({ y_1\plus ty_1}-\x) & S_{122}({ y_1\plus ty_1}-\x)
                     &S_{222}({ y_1\plus ty_1}-\x)  \\
S_{022}({ y_2\plus ty_2}-\x) & S_{122}({ y_2\plus ty_2}-\x)
                       &S_{222}({ y_2\plus ty_2}-\x)  \\
S_{022}({y_3\plus ty_3}-\x) & S_{122}({y_3\plus ty_3}-\x)
                         &S_{222}({y_3\plus ty_3}-\x)  \\
\end{vmatrix}   \, .
\end{equation*}
We have given another expression in \cite{SLC}, separating the variables
$\x$ and $\y$.

The determinant 
$\det\left(
 (tx_i-y_j\frac{1}{t})^{-1} (t^2x_i-y_j\frac{1}{t^2})^{-1}
\right)$
specializes into \\
$\det\left( (x_i-y_j) (x_i^3-y_j^3) ^{-1} \right)$
when $t=\exp(2\pi\sqrt{\moins 1}/3)$. In that case, 
$F_n^1(\x,\y)$ becomes the Schur function in the union of $\x$ and $\y$
of index $[0,0,1,1,2,2,\ldots, n\moins 1,n\moins 1]$ (cf. 
\cite{Stroganov02,Stroganov04}).  
More generally, the Gaudin function of level $r$, when $r$ is odd,
displays such a global symmetry. In that case,       
$(tx-yt^{-1})^{-1}\cdots (t^{r+1}x -yt^{-r-1})^{-1}
                       =  (x-y) (x^{r+2}-y^{r+2})^{-1}$,
and the determinant
$\det\left( (x_i-y_j)(x_i^{r+2}-y_j^{r+2})^{-1} 
\right)$ 
is equal to 
$$  \frac{\Delta(\x) \Delta(\y)}{\prod x_i^{r+2}-y_j^{r+2}}\, 
 S_{0,0,\beta,\beta,\ldots, (n\moins 1)\beta, (n\moins 1)\beta}(\x+\y) \, .$$
More general determinants displaying a symmetry in 
$x_1,\ldots, x_n,y_1,\ldots, y_n$ are given in 
\cite[Lemma 13]{Pfaff} and \cite{Okada}.

We shall now characterize the Gaudin function by specialization.
Expanding the determinant expressing $F^r_n(\x, \y)$ by linearity in $\x$,
one sees that $F^r_n(\x, \y)$ is a linear combination of 
terms $ \frac{1}{\Delta(\y)} 
 \Bigl| y_i^{v_1},\, \ldots,\ldots y_i^{v_n}\Bigr|$, \\  
  $0\leq v_1,\ldots,v_n\leq (n\moins 1)(r\plus 1)$, 
i.e. is a linear combination of Schur functions of $\y$ indexed
by partitions contained in 
  $\boxplus = [\, \underbrace{(n\moins 1)r,\ldots, (n\moins 1)r}_n \, ]$.
By symmetry $\x \leftrightarrow \y$, $F^r_n(\x, \y)$ is 
a linear combination of products of Schur function of $\x$ 
and of Schur functions of $\y$ indexed by partitions 
contained in $\boxplus$.

Given any infinite set of indeterminates $\z$, any linear combination 
of  Schur functions  in $\x$ with coefficients in $\z$, indexed by partitions
$\subseteq \boxplus$, is also a linear combination of 
Grassmannian Schubert polynomials $Y_v(\x,\z)$, $v \subseteq \boxplus$.
As such, it is determined by the 
$\binom{nr+n-r}{n}$ specializations $\x \subset \{z_1,\ldots, z_{nr+n-r}\}$.
In the next theorem, we choose 
$\z=\{ y_1,\ldots, y_n, ty_1,\ldots, ty_n, t^2 y_1,\ldots, t^2 y_n,
\ldots \}$ to get simple specializations.

\begin{theorem}    \label{TheoremSpec}
 $F^r_n(\x,\y)$ is the only linear combination of Schur functions in $\x$,
with coefficients in $\y$, indexed by partitions contained in $\boxplus$,
which has the same specializations 
$$ \x \subset \{ y_1,\ldots, y_n,\ldots, t^r y_1,\ldots, t^r y_n\} $$
than the function
$$ G_n^r(\x,\y) := \frac{\Delta(\x)}{\Delta(\y)} \prod_i  
  S_{\square} (y_i(1\plus \cdots\plus t^r) -\x)  \, ,$$
where $\square$ is, as before, equal to $(n\moins 1)^r$.
\end{theorem} 

\Proof If the specialization of $\x$ contains several occurrences of the same 
$y_i$ (ignoring the powers of $t$), then 
all the functions $S_{j,\, \square}(y_i(1\plus \cdots\plus t^r)-\x)$
vanish. Thus, $F_n^r(\x,\y)$ as well as $G_n^r(\x,\y) $
vanish in that case.

Let now $\x'=\{ x'_1=y_1 t^{\epsilon_1}, \ldots, 
                                 x'_n=y_n t^{\epsilon_n} \}$,
with $0\leq \epsilon_1,\ldots, \epsilon_n \leq r$.
In that case, each $y_i(1\plus \cdots\plus t^r)- \x'$ is equal to
the difference of two sets of respective cardinalities $r,n-1$,
and, according to (\ref{Factorise}),  
$$ S_{j,\, \square}(y_i(1\plus \cdots\plus t^r)-\x')
= S_j( x'_i-\x')\, S_{\square}(y_i(1\plus \cdots\plus t^r)-\x') \, .$$
This factorization allows to extract from the determinant expressing
$F_n^r(\x',\y)$   the factor 
$\prod_i S_{\square} (y_i(1\plus \cdots\plus t^r) -\x')$.
There remains $\det\bigl(S_j( x'_i-\x')  \bigr)$,
which is equal to the Vandermonde $\Delta(\x')$. 
Therefore, the two functions $F_n^r(\x,\y)$,  $G_n^r(\x,\y)$ have the
same specializations in $\{ y_1,\ldots, t^r y_n\} $.
To characterize $F_n^r(\x,\y)$, we need only specialize $\x$ to a subset
of the first $nr+n-r$ elements of $\{ y_1,\ldots, t^r y_n\} $,
and this finishes the proof of the theorem.  \hfill  QED

Notice that all the specializations occurring in Theorem \ref{TheoremSpec} 
are either $0$, or products of factors $(y_it^k -x_j)$.

When $r=1$, 
Theorem \ref{TheoremSpec} claims that $F_n^1(\x,\y)$ has the same 
specializations $\x\subset \{ y_1,\ldots,y_n,ty_1,\ldots, ty_n\}$ 
as $ \Delta(\x) \Delta(\y)^{-1} 
  \prod_i S_{n-1}(y_i+ty_i-\x)$. We have moreover remarked that one can 
suppress one letter from $\{ y_1,\ldots,ty_n\}$
to characterize the Gaudin function. 

For $n=2$, $r=2$, the expression in Theorem \ref{GaudinMatLevel}
specializes, for $x_1=y_1,x_2=t^2y_2$, into
\begin{multline*}
 \frac{1}{y_2-y_1} 
\begin{vmatrix}
S_{011}(ty_1+t^2y_1-t^2y_2)  & S_{111}(ty_1+t^2y_1-t^2y_2) \\
S_{011}(y_2+ty_2-y_1)  & S_{111}(y_2+ty_2-y_1)
\end{vmatrix}  \\
= t^3 (y_2-y_1)(y_1-ty_2)^2 
 \begin{vmatrix}  1 & -t^2y_2 \\   1  & -y_1 
\end{vmatrix}  = t^3 (y_2-y_1)(y_1-ty_2)^2 (t^2y_2-y_1)  \, .
\end{multline*}
Theorem \ref{TheoremSpec} gives on the other hand 
\begin{multline*}
 \frac{t^2y_2-y_1}{y_2-y_1}\, S_{011}( ty_1+t^2y_1-t^2y_2)\, 
                                S_{011}(y_2+ty_2-y_1)      \\
= \frac{t^2y_2-y_1}{y_2-y_1} (ty_1-t^2y_2)(t^2y_1-t^2y_2) (y_2-y_1)(ty_2-y_1)
\, , 
\end{multline*}
which is, indeed, equal.

\section{Euler-Poincar\'e characteristics}

We go back to the original Gaudin-Izergin-Korepin determinant
\cite{Gaudin,Izergin},
that is, from now on, we take level $r=1$.

The \emph{Euler-Poincar\'e characteristics} for a flag manifold
under $Gl(\C^n)$, conveniently generalized by Hirzebruch 
\cite{Hirzebruch, Chiy}, can be combinatorially interpreted as
a summation over the symmetric group $\mfS_n$ ~: 
$$  \C[t][x_1,\ldots,x_n] \ni f \longrightarrow 
  f \carre_\omega := \sum_{w\in\mfS_n} \left( 
  f \frac{\prod_{i<j} (tx_i -x_j)}{x_i-x_j}  \right)^w  \in \Sym(\x)  \, . $$

This morphism is characterized by the fact that the images of
dominant monomials 
$$   x^\l := x_1^{\l_1} \cdots x_n^{\l_n}\, ,\ \l_1\geq \cdots \geq \l_n\geq 0
$$
are the Hall-Littlewood polynomials \cite{Macdonald}
 $c_\l P_\l(\x,t)$, the normalization
constants $c_\l$, writing $\l= 0^{m_0} 1^{m_1}\cdots n^{m_n}$, being 
$$   c_\l = \prod_{i=0}^n \prod_{j=1}^{m_i} (1-t^j)(1-t)^{-1} \, .  $$

The elementary operators (case $n=2$) are 
$$ \carre_i := \carre_{s_i} :
  f \longrightarrow f\, \carre_i = f (tx_i-x_{i+1})\, \d_i $$
and generate the Hecke algebra of the symmetric group (as an algebra
of operators on polynomials. This is the description
of the affine Hecke algebra that I used with M.P. Sch\"utzenberger
in \cite{SymmetryFlag}).
The usual generators  
$  T_i := \carre_i -1  $
satisfy the braid relations, and the Hecke relation $(T_i-t)(T_i+1)=0$,
while $\carre_i^2 = (1+t)\carre_i$.

We shall also need an affine operation $\theta$, which is
the incrementation of indices on the $x$-variables :
\begin{equation*}  
x_i\, \theta =x_{i+1}\, , \quad \text{periodicity}\ x_{i+n} = x_i t^{-1} \, .
\end{equation*}
Notice that to define Macdonald's polynomials (see next section), 
one uses the periodicity
$x_{i+n} = q x_i$ with $q$ independent of $t$.

\begin{theorem}
Let $\x,\y$ be two alphabets of cardinality $n$, $f$ be a function
of a single variable. Then
\begin{multline}     \label{TheoremTheta}
 f(x_1)\, R(\x-x_1\, ,\,  \y) (1-t\theta)\cdots (1-t^{n-1}\theta) 
                                   \carre_\omega  \\
 = \left(f(x_1) x_2\cdots x_n \d_1\cdots \d_{n-1}\,\right)
                                                F_n^1(\x,\y)\, [n]! \, .
\end{multline}

\end{theorem}

\Proof The LHS, as a function of $\y$, belongs to the space
generated by the Schur functions of index contained in the 
partition $(n\moins 1)^n$. It therefore can be determined
by computing all the specializations 
$$ \y \subset \{ x_1,\ldots, x_n, x_2 t^{-1},\ldots  x_n t^{-1} \} $$
(we do not need to take $x_1t^{-1}$).

To lighten notations, let us take $n=4$.  The function 
$$ f(x_1)\, R(x_2\plus x_3\plus x_4 ,\y) (1-t\theta)
  (1-t^2\theta)(1-t^3\theta) $$
is equal to 
\begin{multline*}
 f(x_1)\, R(x_2\plus x_3\plus x_4 ,\y)
 -t {3 \brack 1} f(x_2)\, R(x_3\plus x_4\plus x_1/t ,\y)   \\
 +t^3 {3 \brack 2} f(x_3)\, R(x_4\plus x_1/t\plus x_2/t ,\y)
 -t^6 f(x_4)\, R(x_1/t\plus x_2/t\plus x_3/t ,\y)
\end{multline*}

The sum under the symmetric group can be written
\begin{multline*}
 \sum_w \left(f(x_1)R(x_2\plus x_3\plus x_4 ,\y) 
                      \frac{\Delta_t(1234)}{\Delta(1234)}   \right)^w \\ 
 -t {3 \brack 1} \sum_w \left(f(x_1)R(x_3\plus x_4\plus x_2t^{-1} ,\y)
                      \frac{\Delta_t(2134)}{\Delta(2134)}   \right)^w  \\
 +t^3 {3 \brack 2} \sum_w \left(f(x_1)R(x_4\plus x_2 t^{-1}\plus x_3t^{-1} ,\y)
                      \frac{\Delta_t(3214)}{\Delta(3214)}   \right)^w \\
 -t^6 \sum_w \left(f(x_1)R( x_2 t^{-1}\plus x_3t^{-1}\plus x_4t^{-1} ,\y)
                      \frac{\Delta_t(4231)}{\Delta(4231)}   \right)^w
\end{multline*}

We shall identify the coefficients of $f(x_1)$ in both members of 
\ref{TheoremTheta}. 
Consider all the specializations 
$\y \subset \{ x_1,\ldots, x_4, x_2/t,\ldots, x_4/t \}   $
of the LHS.  
Up to symmetry, the only non-zero specializations are
\begin{itemize}
\item  $\y\to \{x_1,x_2,x_3,x_4\}$:\, 
  $-t^6 R(\frac{x_2+x_3+x_4}{t},\y)\sum_{w:\, w_1=1} 
            \left( \frac{\Delta_t(1234)}{\Delta(1234)}  \right)^w $
\item  $\y\to \{x_1,x_2,x_3,\frac{x_4}{t} \}$:\,    
 $ t^3 {3 \brack 2} R(\frac{x_2+x_3}{t}+x_4,\y)  \sum_{w:\, w_1=1,w_4=4}
  \left( \frac{\Delta_t(3214)}{\Delta(3214)}   \right)^w   $   
\item $\y\to \{x_1,x_2,\frac{x_3}{t},\frac{x_4}{t}   \}$:\,
 $-t {3 \brack 1} R(\frac{x_2}{t}+x_3\plus x_4,\y)  \sum_{w:\, w_1=1,w_2=2}
  \left( \frac{\Delta_t(2134)}{\Delta(2134)}   \right)^w $
\item $\y\to \{x_1,\frac{x_2}{t},\frac{x_2}{t},\frac{x_4}{t}    \}$:\,
 $ R(x_2\plus x_3\plus x_4,\y)  \sum_{w:\, w_1=1}
   \left( \frac{\Delta_t(1234)}{\Delta(1234)}   \right)^w $\, .
\end{itemize} 
Up to the global factor
$$ \frac{x_2x_3x_4(1-t)(1-t^2)(1-t^3)}{\Delta_t(1234) \Delta_t(4321)
 R(x_1, x_2\plus x_3\plus x_4) }    $$
these specializations are respectively equal to
$$ 1,\, \frac{R(x_1\plus x_2\plus x_3, x_4)}{R(x_1\plus x_2\plus x_3, tx_4)},\, 
\frac{R(x_1\plus x_2,x_3\plus x_4)}{R(x_1\plus x_2,tx_3\plus tx_4)},\,
\frac{R(x_1, x_2\plus x_3\plus x_4)}{R(x_1, tx_2\plus tx_3\plus tx_4)}\, . $$
They coincide with the specializations of the RHS of (\ref{TheoremTheta}),
thanks to Theorem \ref{TheoremSpec}, 
writing $f(x_1)x_2x_3x_4 \d_1\d_2\d_3$ as 
$$   \frac{f(x_1) x_2x_3x_4}{R(x_1, x_2\plus x_3\plus x_4)}
 + \frac{f(x_2) x_1x_3x_4}{R(x_2, x_1\plus x_3\plus x_4)}
+ \frac{f(x_3) x_1x_2x_4}{R(x_3, x_1\plus x_2\plus x_4)} 
+ \frac{f(x_4) x_1x_2x_3}{R(x_4, x_1\plus x_2\plus x_3)}  \, . $$

In final, we have checked enough specializations to prove 
(\ref{TheoremTheta}).      \QED

\section{Generating functions of Macdonald polynomials}

The symmetric Macdonald polynomials $P_\l(\x;q,t)$
satisfy a Cauchy formula~:
\begin{equation}  \label{CauchyMacdo}
   \sigma_1\left(\x\y\frac{1-t}{1-q} \right)  :=
  \prod_{x\in\x, y\in\y} \prod_{i\geq 0} 
     \frac{1-tq^i xy}{1-q^i xy} = 
               \sum_\l b_\l P_\l(\x;q,t) P_\l(\y;q,t) \, ,
\end{equation}
sum over all partitions of length $\ell(\l)\leq n$,
the constants $b_\l $ being defined in \cite[VI.4.11]{Macdonald}.

Let $\tau_q$ be the following incrementation of indices on the $x$-variables :
\begin{equation*}
x_i\, \tau_q =x_{i+1}\, , \quad \text{periodicity}\ x_{i+n} = qx_i  \, .
\end{equation*}

We want to compute 
$$ \sigma_1\left(\x\y\frac{1-t}{1-q} \right)\, 
   (1-t\tau_q)\cdots (1-t^n\tau_q)\, \carre_\omega \, . $$
Since 
$$ \x \frac{1-t}{1-q} \tau_q = \x \frac{1-t}{1-q} + x_1(t-1)  \, , $$
one has
\begin{multline}   \label{fgMacdo}
 \sigma_1\left(\x\y\frac{1-t}{1-q} \right)\,
   (1-t\tau_q)\cdots (1-t^n\tau_q)\, \carre_\omega   \\
 = \sigma_1\left(\x\y(1\moins t)\right) (1-t\tau_0)\cdots (1-t^n\tau_0)
           \carre_\omega \,  \sigma_1\left(\x\y q\frac{1-t}{1-q} \right)
\end{multline}

The parameter $q$ has been eliminated from the operation, 
and we are thus reduced to 
the case of Hall-Littlewood polynomials, which is treated in the next
theorem.

\begin{theorem}   \label{ThHL}
The image of the generating function of Hall-Littlewood polynomials 
$P_\l(\x,t)= P_\l(\x;0,t)$ under 
$(1-t\tau_0)\cdots (1-t^n\tau_0) \carre_\omega$ is 
\begin{equation}                                        \label{fgHL}
\sigma_1\left(\x\y(1\moins t)\right) (1-t\tau_0)\cdots (1-t^n\tau_0)
           \carre_\omega 
 =  \sigma_1(\x\y)\,  \tF_n^1(\x,\y)\, [n]!  \, ,
\end{equation}
where $\tF_n^1(\x,\y)$ is the Gaudin function
$(x_1\cdots x_n)^{n-1}\, F_n^1(\x^\vee, \y)$, \\  
$\x^\vee=\{ x_1^{-1},\ldots, x_n^{-1}\}$, and $[n]!=(1-t)\cdots (1-t^n)$. 
\end{theorem}

\Proof
One rewrites 
$ \sigma_1\left(\x\y(1\moins t)\right) 
                   = R(t\x, \y^\vee) R(\x, \y^\vee)^{-1}$. 
Notice that 
\begin{multline*}
  \frac{ R(tx_i+\cdots+tx_n, \y^\vee)}{ R(x_i+\cdots+x_n, \y^\vee)}  
  (1-t^n\tau_0) = 
\frac{ R(tx_i+\cdots+tx_n, \y^\vee)}{ R(x_i+\cdots+x_n, \y^\vee)}
-\frac{ R(tx_{i+1}+\cdots+tx_n, \y^\vee)}{ R(x_{i+1}+\cdots+x_n, \y^\vee)}\\
= \bigl( R(tx_i,\y^\vee)-t^n R(x_i,\y^\vee) \bigr) 
 \frac{ R(tx_{i+1}+\cdots+tx_n, \y^\vee)}{ R(x_i+\cdots+x_n, \y^\vee)}
= \frac{f(x_i)}{R(x_i+\cdots+x_n, \y^\vee)} \, 
\end{multline*}
with $f(x_i)$ a polynomial in $x_i$ of degree $n-1$.
Multiplying the LHS of (\ref{fgHL}) by the function 
$R(\x,\y^\vee)$, one transforms it into
\begin{multline*}
\Bigl( f(x_1) R(tx_2+\cdots +tx_n,\y^\vee) -e_1 f(x_2) R(tx_3 
                                           +\cdots +tx_n+x_1,\y^\vee) \\
 +\cdots +(\moins 1)^{n-1} e_{n-1} 
            f(x_n) R(x_1+\cdots +x_{n-1}, \y^\vee) \Bigr)\, \carre_\omega\, ,
\end{multline*}
where $e_1,\ldots, e_{n-1}$ are the elementary symmetric functions
of $t,\ldots, t^{n-1}$. 
One recognizes in this last expression 
$$  f(x_1) R(tx_2+\cdots +tx_n,\y^\vee) (1-t\theta)\cdots 
(1-t^{n-1}\theta) \, . $$

We now invoke Theorem \ref{TheoremTheta}.
Since the function $f(x_1) x_2\cdots x_n \d_1\cdots \d_{n-1}$ is a constant 
(for degree reasons), the LHS of (\ref{fgHL}) is 
proportional to $F_n^1(\x,\y^\vee)$, that is, is
proportional to $\tF_n^1(\x,\y)$. 
In fact, $f(x) = \prod_j( tx -y_j^{-1})\\  -t^n \prod_j( x -y_j^{-1})$,
so that 
$$ f(x_1) x_2\cdots x_n \d_1\cdots \d_{n-1} 
                       = (1-t^n)(y_1\cdots y_n)^{-1} \, .  $$

Correcting by the right powers of $(x_1\cdots x_n)$ and 
$(y_1\cdots y_n)$, one finishes the proof of the theorem.  \QED

\bigskip
One can now go back to the case of Macdonald polynomials, and recover
a result of Warnaar \cite[Th.3.1]{Warnaar}.

\begin{theorem}[Warnaar]
There holds 
\begin{equation}  \label{Warnaar}
 \sum_\l b_\l P_\l(\x;q,t) P_\l(\y;q,t) 
  \prod_{i=1}^n (1-q^{\l_i} t^{n-i+1})  = 
   \sigma_1\left(\x\y\frac{1-t}{1-q} \right)\, 
  \sigma_1(t\x\y) \tF_n^1(\x,\y) \, .
\end{equation}
\end{theorem}

\Proof
The non symmetric Macdonald polynomials are eigenfunctions
of certain commuting Dunkl-type operators $\xi_1,\ldots, \xi_n$,
  first introduced in \cite{BGHP} and extensively used by 
Cherednik \cite{Cherednik}.
The eigenvalues are $q^{\l_1} t^{n-1},\ldots, q^{\l_n}t^0$
for the polynomial $M_\l$ indexed by $\l: \l_1\geq \cdots \geq \l_n\geq0$.
Up to normalization, the image of $M_\l$ under $\carre_\omega$ is
equal to $P_\l(\x,q,t)$, and for any symmetric function $g$ in
$n$ variables, then  
$$ P_\l(\x;q,t)\, g(\xi_1,\ldots, \xi_n)
      =   P_\l(\x;q,t)\, g( q^{\l_1} t^{n-1},\ldots, q^{\l_n}t^0)  \, . $$

Using \cite{YangRota} that $\carre_\omega \, \xi_i\, \carre_\omega =
 \carre_\omega \, t^{n-i}\tau_q\, \carre_\omega$,
one sees that 
$$  P_\l(\x;q,t) \prod_{i=1}^n (1-q^{\l_i} t^{n-i+1})  = 
   P_\l(\x;q,t) (1-t\tau_q)\cdots (1-t^n\tau_q)\, \carre_\omega \, .$$ 
Therefore, the LHS of (\ref{Warnaar}) can be identified with
$$\sigma_1\left(\x\y \frac{1\moins t}{1\moins q} \right) 
            (1-t\tau_q)\cdots (1-t^n\tau_q)\, \carre_\omega    \, .$$
Thanks to (\ref{fgMacdo}) and (\ref{fgHL}), 
 this can be written 
$$  \sigma_1(\x\y) \sigma_1\left(\x\y q\frac{1\moins t}{1\moins q}\right) 
    \tF_n^1(\x,\y)
= \sigma_1\left(\x\y \frac{1\moins t}{1\moins q}  \right) 
             \sigma_1(t\x\y) \tF_n^1(\x,\y) \, ,$$
which is Warnaar's formula.                    \QED

\section{Note: Symmetric functions and Schubert polynomials}

We use $\lambda$-ring conventions to describe symmetric functions.
Given three sets $A,B,C$ of indeterminates (``alphabets''),
the generating function of complete functions $S_n(AB-C)$ is 
$$\sigma_z:= \frac{\prod_{c\in \C} 1-zc}{\prod_{a\in A,b\in B} 1-z ab }
= \sum z^n S_n(AB-C)\, .$$
We write alphabets as sums of the letters composing them. 
For example, $A(1+t+\cdots+t^r)$ is the alphabet $\{ at^i :\, a\in A,
0\leq i\leq r \}$. 

Schur functions $S_v(A-C)$, $v\in \N^n$, are determinants of complete
functions:
$$   S_v(A-C) = \det\bigl( S_{v_j+j-i}(A-C)\bigr)_{i,j=1\ldots n}  \, .$$

One generalizes Schur functions to \emph{multi-Schur functions} by taking
 different alphabets in blocks of columns of the preceding determinant.
For example, $S_{v_1; v_2,\ldots, v_n}(A_1-C_1; A-C)$ 
is the determinant with first column
$S_{v_1}(A_1-C_1),\ldots, S_{v_1-n+1}(A_1-C_1)$,
 and entries $S_{v_j+j-i}(A-C)$ elsewhere. 

Multi-Schur functions satisfy some factorization properties
\cite[Prop. 1.4.3]{Cbms}. We need only the following case, which was much used
in classical elimination theory in the 19th century.

\begin{lemma}
Given two finite alphabets $A,B$ of respective cardinalities $\alpha,\beta$,
given $j:\, 0\leq j\leq \beta$, then 
\begin{equation}  \label{Factorise}
  S_{j,\beta^{\alpha}}(A-B)=  S_j(-B) \prod_{a\in A, b\in B} (a-b)=
 (\moins 1)^j e_j(B) \prod_{a\in A, b\in B} (a-b)  \, ,
\end{equation}
where $e_j(B)$ is the elementary symmetric function of degree $j$ in $B$.
In particular, $S_{j,\beta^{\alpha}}(A-B)=0$ if $A\cap B\neq \emptyset$.
\end{lemma}

There are several families of non-symmetric polynomials extending the
basis of Schur functions.
Of special interest are the \emph{Schubert polynomials}
$Y_v(\x,\y)$, $v\in \N^n$, which constitute a linear basis of the
ring of polynomials in $x=\{x_1,\ldots, x_n\}$, with coefficients in
$y_1,y_2,\ldots, y_\infty$. 
They can be characterized by vanishing properties related to the Bruhat order.

The subfamily of Schubert polynomials indexed by (increasing) partitions
form a basis of the ring of symmetric polynomials \cite{Cbms}. 
It satisfies the following property.
Given $u\in \N^n:\, 0\leq u_1 \leq \cdots \leq u_n$, let
$\y^{<u>} :=\{ y_{1+u_1},\ldots, y_{n+u_n}\}$. 

\begin{lemma}
Let $v= [0\leq v_1\leq \cdots \leq v_n]$ be a partition.
 Then $Y_v(\x,\y)$ is the only 
symmetric function in $\x$ of degree $|v|= v_1+\cdots+v_n$ such that
$ Y_v(\y^{<u>},\y)=0$ for all $u:\, |u|\leq |v|$, $u\neq v$,
and $Y_v(\x, \{0,0,\ldots\}) = S_v(\x)$.  

More generally, $ Y_v(\y^{<u>},\y) \neq 0$ iff the diagram 
of $u$ contains the diagram of $v$.
\end{lemma}

In the case where $\y=\{0,1,2,\ldots\}$ (resp.
$\y=\{q^0,q^1,q^2,\ldots\}$), the polynomials
$Y_v(\x,\y)$ are called \emph{factorial Schur functions}
(resp. \emph{$q$-factorial Schur functions}, and the above
vanishing properties are extensively used in \cite{OkounkovOlshanski}.

\subsection*{Acknowledgements}
The author benefits from the ANR project BLAN06-2\_134516.

\bigskip
\begin{center}
{\large Alain Lascoux}\\
 CNRS, IGM, Universit\'e Paris-Est \\
  77454 Marne-la-Vall\'ee Cedex, France\\
  Email: Alain.Lascoux@univ-mlv.fr\\
\end{center}

\end{document}